\documentclass{amsart}

\usepackage[mathcal]{eucal}

\usepackage{amssymb, graphics}

\usepackage[all]{xy}

\newtheorem{thm}{Theorem}

\newtheorem{lem}[thm]{Lemma}
\newtheorem{prop}[thm]{Proposition}
\theoremstyle{definition}

\theoremstyle{remark}
\newtheorem{rem}[thm]{Remark}
\numberwithin{equation}{section}

\begin{document}

\title[Quantum hyperbolic invariants]
{Quantum hyperbolic invariants \\
for diffeomorphisms of small surfaces}

\author{Xiaobo Liu}

\address{Department of Mathematics,
Columbia University, New York, NY 10027, U.S.A.}

\email{xiaoboli@math.columbia.edu}

\thanks{ This work was partially supported by NSF
grant DMS-0103511 at the University of Southern California.}

\date{\today}


\begin{abstract}
An earlier article \cite{BonLiu} introduced new invariants for
pseudo-Anosov diffeomorphisms of surface, based on the
representation theory of the quantum Teichm\"uller space. We
explicity compute these quantum hyperbolic invariants in the case of
the 1--puncture torus and the 4--puncture sphere.
\end{abstract}

\maketitle

\section{Introduction}

In \cite{BonLiu}, Francis Bonahon and the author constructed
quantum hyperbolic invariants for pseudo-Anosov diffeomorphisms
of a punctured surface
$S$. For every odd integer $N$, these invariants associate to the
pseudo-Anosov diffeomorphism
$\varphi: S\to S$ a square matrix $C_\varphi$ of dimension
$N^{3g+p-3}$ (where $g\geq0$ is the genus of $S$ and $p>0$ is its
number of punctures), defined up to conjugation and scalar
multiplication. In particular, $C_\varphi$ is an $N\times N$
matrix for the 1--puncture torus and the 4--puncture sphere. The
current paper provides an explicit computation for these two
surfaces. It is based on the fact that the mapping class group of
these surfaces is particularly simple.

For the 1-puncture torus, the isotopy class of a diffeomorphism
$\varphi: S\to S$ is completely determined by the matrix
$A_\varphi\in\mathrm{SL}_2(\mathbb Z)$ defined by considering the
action of $\varphi$ on $H_1(S)=\mathbb Z^2$. A diffeomorphism
$\varphi$ is pseudo-Anosov if and only if $\left\vert
\mathrm{Tr}(A_\varphi)\right\vert >2$. Such a matrix $A_\varphi$ is
conjugate to a product
$$A_\varphi=A_1A_2\cdots A_n,\quad\textrm{where}\quad
A_i=R=\left(\begin{array}{cc} 1 & 1\\0 & 1\end{array}\right)\
\textrm{or}\ A_i=L=\left(\begin{array}{cc} 1 & 0\\1 &
1\end{array}\right).$$

Our computation makes use of the matrices $C_R(u,v,u',v',h)$ and
$C_L(u,v,u',v',h)$ defined by
\begin{align*} {} & C_R(u,v,u',v',h)_{ij} \\
= & \ q^{2(j-i)^2}\ \Big(\frac{u'v'}{uh}\Big)^{j-i}\
\Big(\frac{v'}{v}\Big)^i\
\prod_{\alpha=1}^{i}\frac{1}{\big(1+q^{4\alpha-3}u\big)
\big(1+q^{4\alpha-1}u\big)},
\end{align*}
and
\begin{align*}
{} & C_L(u,v,u'',v'',h)_{ij} \\
= & \ \sum_{k=1}^N q^{4ik}\ q^{2(j-k)^2+2k^2}\
\Big(\frac{u''v''}{vh}\Big)^{j-k}\ \Big(\frac{uvv''}{h}\Big)^k\
\prod_{\alpha=1}^{k}\frac{1}{\big(1+q^{4\alpha-3}v\big)
\big(1+q^{4\alpha-1}v\big)}.
\end{align*}

\begin{thm}
Let $\varphi : S \to S$ be a pseudo-Anosov diffeomorphism of the
$1$--puncture torus $S$, and consider a matrix
$A_\varphi=A_1A_2\cdots A_n$, with
$A_i=R$ or
$L$, associated to $\varphi$ as above. Then
$$C_\varphi=C_1 C_2\cdots C_n $$ where
$C_i=C_R(u_{i-1},v_{i-1},u_i,v_i,1)$ if $A_i=R$, $
C_i=C_L(u_{i-1},v_{i-1},u_i,v_i,1)$ if $A_i=L$, and where the
complex numbers $u_i$, $v_i$ are explicitly defined in
\S\ref{sec:AlgLem} and \S\ref{sec:1-PuncTorus} and are determined by
the complete hyperbolic metric of the mapping torus $M_\varphi$ of
$\varphi$.
\end{thm}
There is a similar theorem in the case of the 4-puncture sphere,
using different matrix functions $C^*_R$ and $C^*_L$ defined in
\S\ref{sec:4-PuncSphere}. In this situation, the action on slopes
defines a map from the mapping class group of the 4-puncture sphere
to $\mathrm{PSL}_2(\mathbb Z)$ with kernel $\mathbb Z/2\mathbb Z\
\oplus\ \mathbb Z/2\mathbb Z$. The image of a diffeomorphism
$\varphi$ can be represented by a matrix
$A_\varphi\in\mathrm{SL}_2(\mathbb Z)$. Again, $\varphi$ is
pseudo-Anosov if and only if $|\mathrm{Tr}(A_\varphi)|>2$.
\begin{thm}
Let $\varphi : S \to S$ be a pseudo-Anosov diffeomorphism of the
$4$--puncture sphere $S$, and consider a matrix
$A_\varphi=A_1A_2\cdots A_n$, with
$A_i=R$ or
$L$, associated to $\varphi$ as above. Then
$$C_\varphi=C^*_1 C^*_2\cdots C^*_n $$ where
$C^*_i=C^*_R(u_{i-1},v_{i-1},u_i,v_i,1)$ if $A_i=R$, $
C^*_i=C^*_L(u_{i-1},v_{i-1},u_i,v_i,1)$ if $A_i=L$, and where the
complex numbers $u_i$, $v_i$ are explicitly defined in
\S\ref{sec:4-PuncSphere} and are determined by the complete
hyperbolic metric of the mapping torus $M_\varphi$ of $\varphi$.
\end{thm}

\noindent\textbf{Acknowledgements:} We would like to thank
Fran\c cois Gu\'eritaud for crucial help with the combinatorics
of the 1--puncture torus, and Francis Bonahon for general
advice and encouragement.

\section{Quantum hyperbolic invariants of surface diffeomorphisms}

We briefly sketch the construction in \cite{BonLiu} of
quantum hyperbolic invariants for pseudo-Anosov diffeomorphisms
of a punctured surface. Details will be provided in later
sections in the specific case of the 1--puncture torus and the
4--puncture sphere.

Let $S$ be a punctured surface, that
is, a closed surface with finitely many  points removed. The
Teichm\"uller space of
$S$, denoted by $\mathcal T_S$, is the space of isotopy classes
of complete hyperbolic metrics on $S$. An ideal triangulation
$\lambda$ provides a certain set of global coordinates for
$\mathcal T_S$, called the exponential shear coordinates
(these essentially are cross-ratios), which parameterize this
space by a convex cell
$\mathbb R_+^n$ (more accurately they parametrize a
finite-sheeted branched covering of $\mathcal T_S$ called the
enhanced Teichm\"uller space).

The \emph{quantum Teichm\"uller space} of $S$ is a noncommutative
deformation of the algebra of rational functions on $\mathcal
T_S$, with a deformation parameter $q\in\mathbb C$; see
\cite{FocChe, Kas, BonLiu, Liu}. The construction is based on a
non-commutative analog of the shear coordinate parametrization
for the quantum Teichm\"uller space, the
\emph{Chekhov-Fock algebra}
$\mathcal T^q_\lambda$ associated to an ideal triangulation
$\lambda$ of $S$.

When we change the ideal
triangulation $\lambda$ to another one
$\lambda'$, the exponential shear coordinates transform
rationally. There is a non-commutative procedure that
mimics  this transformation rule of cross-ratios: it consists
of an algebraic isomorphism $\Phi^q_{\lambda\lambda'}:
\widehat{\mathcal T}^q_{\lambda'} \to \widehat{\mathcal
T}^q_\lambda$, where
$\widehat{\mathcal T}^q_\lambda$ is the division algebra of formal
fractions of $\mathcal T^q_\lambda$. This \emph{quantum
coordinate change isomorphism} was originally proposed by Chekhov
and Fock, and expressed in terms of a holomorphic function called
the quantum dilogarithm (see \cite{FocChe}).

The Chekhov-Fock algebra $\mathcal T^q_\lambda$ admits finite
dimensional representations only when $q$ is a root of unity.
Assume that
$q$ is a primitive $N$--root of unity with $N$ odd. Any finite
dimensional irreducible representation of this algebra has
dimension $N^{3g+p-3}$ (where $g\geq0$ is the genus of $S$ and $p>0$ is its
number of punctures) and, up to finitely many choices,
is determined by complex weights $x_e\in \mathbb C-\{0\}$
associated to the edges $e$ of the ideal triangulation
$\lambda$.

It turns out that this representation theory of the Chekhov-Fock
algebras $\widehat{\mathcal
T}^q_\lambda$ is well-behaved with respect to the coordinate
change isomorphisms $\Phi^q_{\lambda\lambda'}$. Namely, suppose
that the irreducible
representation
$\rho: \mathcal T^q_\lambda\to \mathrm{End}(V)$ is associated
to  weights $x_e$, as $e$ ranges over all edges of $\lambda$,
and that these weights are in the domain of the
(complexification of the) crossratio transformation rule between
the exponential shear coordinates associated to $\lambda$ and
those associated to $\lambda'$. Then it is possible to make
sense of the irreducible representation
$\rho'=\rho\circ\Phi^q_{\lambda\lambda'}: \mathcal
T^q_{\lambda'}\to \mathrm{End}(V)$, and  $\rho'$ is classified
by the edge weights on $\lambda'$ that are the image of the
weights $x_e$ under the cross-ratio transformation rule.

Now, consider a pseudo-Anosov diffeomorphism $\varphi:S\to S$.
Thurston's Hyperbolization Theorem and Mostow's Rigidity Theorem
provide a unique finite-volume complete hyperbolic metric on the
mapping torus $M_{\varphi}=S\times[0,1]/\sim $, where $\sim$
identifies each $(x,0)$ with $(\varphi(x),1)$. An arbitrary ideal
triangulation $\lambda$ determines a unique pleated surface
$f_\lambda:S\to M_\varphi$ with pleating locus $\lambda$, namely
such that $f$ sends each face of $\lambda$ to a totally geodesic
triangle in $M_\varphi$. This pleated surface defines a complex
weight $x_e\in \mathbb C-\{0\}$ for each edge $e$ of $\lambda$,
called its exponential shear-bend parameter. If we identify
$\lambda$ and $\lambda' = \varphi(\lambda)$, these shear-bend
coordinates have the property that they are invariant under the
cross-ratio transformation rule from $\lambda$ to $\lambda'$. Using
the correspondence between edge weights and representations of the
Chekhov-Fock algebra, this gives a representation $\rho: \mathcal
T^q_\lambda\to \mathrm{End}(V)$ such that
$\rho'=\rho\circ\Phi^q_{\lambda\lambda'}$ is isomorphic to
$\rho\circ \Phi$, using the identification $\Phi:\mathcal
T^q_{\lambda'}\cong \mathcal T^q_\lambda$ provided by $\varphi$. In
other words, there exists an isomorphism $C_\varphi\in
\mathrm{GL}(V)$ such that
\begin{equation}\label{eqn:invariant}
\rho\circ\Phi^q_{\lambda\varphi(\lambda)}(X)=C_\varphi\cdot
\Big(\rho\circ\Phi(X)\Big)\cdot C_\varphi^{-1}\qquad \textrm{for
any}\quad X\in\mathcal T^q_{\varphi(\lambda)}.
\end{equation}

\begin{thm}[Theorem 40 in \cite{BonLiu}] \label{thm:invariant}
The isomorphism $C_\varphi\in \mathrm{GL}(V)$ defined above
depends only on $q$ and $\varphi$, up to conjugation and scalar
multiplication. \qed
\end{thm}

This article is devoted to computing the invariant $C_\varphi$ when
the surface $S$ is  the 1--puncture torus or the 4-puncture sphere.
In these cases, the irreduclible representations of the Chekhov-Fock
algebra have dimension $N$, so that $C_\varphi$ can be considered as
an $N\times N$ matrix defined up to conjugation and scalar
multiplication.

\section{Quantum Teichm\"uller space
of the 1--puncture torus} 
\begin{figure}[h]
\begin{center}
\includegraphics{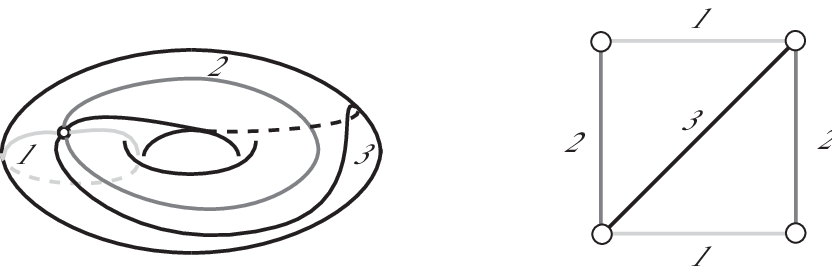}
\end{center}
\caption{An ideal triangulation $\lambda$. In this ideal
triangulation $\sigma_{12}=\sigma_{23}=\sigma_{31}=-2$.}
\label{fig:PTorus}
\end{figure}

Following the construction of the quantum
enhanced Teichm\"uller space \cite{Liu}, fix an \emph{ideal
triangulation} $\lambda$ of the 1-puncture torus $T$, namely a
triangulation of the unpunctured torus $\bar{T}$ with a single
vertex, located at the puncture. The \emph{Chekhov-Fock algebra}
$\mathcal T^q_{\lambda}$ is defined by generators $X_1^{\pm1}$,
$X_2^{\pm1}$, and $X_3^{\pm1}$, respectively associated
  to the edges $\lambda_1$, $\lambda_2$,
$\lambda_3$ of $\lambda$, and by relations determined by the
topology of $\lambda$ as follows:
\begin{equation*}
X_iX_j=q^{2\sigma_{ij}}X_jX_i
\end{equation*}
where $\sigma_{ij}\in\ \{-2,-1,0,1,2\}$
is the number of times one sees \includegraphics{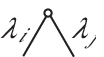} minus
the number of times one sees
\includegraphics{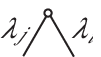}. In particular, the $\sigma_{ij}$
are antisymmetric in the subscripts.
When the $\lambda_i$ occur counterclockwise around the faces of
$\lambda$ as in  Figure~\ref{fig:PTorus}, then
$\sigma_{ij}=-2$ whenever
$j=i+1 \mod 3$.

\begin{figure}[h]
\begin{center}
\includegraphics{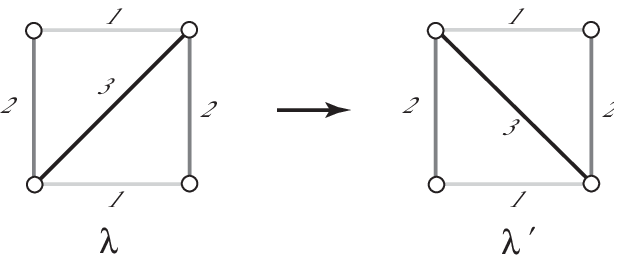}
\end{center}
\caption{The diagonal exchange $\Delta_3$ on edge $\lambda_3$.}
\label{fig:DiaEx}
\end{figure}

The \emph{diagonal exchange} $\Delta_3$ takes the ideal
triangulation $\lambda$ to a new ideal triangulation $\lambda'$
obtained by replacing
$\lambda_3$ by the other diagonal $\lambda'_3$ of the square
formed by $\lambda_1$ and $\lambda_2$, as Figure~\ref{fig:DiaEx}.
If $\widehat{\mathcal{T}}^q_{\lambda}$ denotes the fraction
division algebra of $\mathcal
T^q_\lambda$, consisting of rational functions in the
skew-commuting variable $X_1$, $X_2$, $X_3$, the coordinate
change isomorphism
$\Phi^q_{\lambda\lambda'}:\widehat{\mathcal{T}}^q_{\lambda'}
\to \widehat{\mathcal{T}}^q_{\lambda}$
introduced in \cite{Liu} is such that
\begin{align} \label{eqn:diaex}
\Phi^q_{\lambda\lambda'}(X_1')=&(1+qX_3)(1+q^3X_3)X_1,\nonumber\\
\Phi^q_{\lambda\lambda'}(X_2')=&\big(1+qX_3^{-1} \big)^{-1}
\big(1+q^3X_3^{-1} \big)^{-1}X_2,\\
\Phi^q_{\lambda\lambda'}(X_3')=&X_3^{-1}.\nonumber
\end{align}
in the case of Figure~\ref{fig:PTorus}, namely if
$\sigma_{ij}=-2$ when
$j=i+1 \mod 3$.

Now suppose that $q$ is a primitive $N$-th root of unity with
$N$ odd. The center of $\mathcal T^q_\lambda$ is generated by
$X_1^N$,
$X_2^N$, and $X_1X_2X_3$.

\begin{thm}[A special case of Theorem 21 in \cite{BonLiu}]
\label{thm:repn}
Let $q$ be a primitive $N$-th root of unity with $N$ odd. Any
irreducible representation
$\rho$ of the Chekhov-Fock algebra is of dimension $N$, and the
conjugacy class of $\rho$ is determined by its restriction to the
central elements $X_1^N$, $X_2^N$ and
$H=q^{-(\sigma_{12}+\sigma_{13}+\sigma_{23})}X_1X_2X_3$. \qed
\end{thm}
The coefficient introduced in the definition
$H=q^{-(\sigma_{12}+\sigma_{13}+\sigma_{23})}X_1X_2X_3$ is
designed to make $H$ independent of the order of the
$\lambda_i$. This will also guarantee that $H$
is invariant under coordinate change isomorphisms in the sense
that
$\Phi^q_{\lambda\lambda'}(H')=H$ as proved in \cite{Liu}. Also,
note that $\sigma_{12}+\sigma_{13}+\sigma_{23} =\pm2$, which may
make this exponent less intimidating.

If $\rho:\mathcal T^q_\lambda\to\mathrm{End}(V)$ is an irreducible
representation, then $\rho(X_1^N)=x_1\,\mathrm{id}_V$,
$\rho(X_2^N)=x_2\,\mathrm{id}_V$, and $\rho(H)=h\, \mathrm{id}_V$
for complex numbers $x_1,x_2,h\in\mathbb C-\{0\}$. It follows that
$\rho(X_3^N)=x_3\,\mathrm{id}_V$ where $x_3=h^N/(x_1x_2)$. In the
next section we shall interpret $(x_1,x_2,x_3)$ as the shear-bend
coordinates of a pleated 1-puncture torus. The following lemma
provides a base for this interpretation. Letting $q=1$ in
Formula~(\ref{eqn:diaex}) so that the variables commute, the
isomorphism $\Phi^1_{\lambda\lambda'}$ defined there can be seen as
a rational function from $\big(\mathbb C-\{0,-1\}\big)^3$ to itself.
This $\Phi^1_{\lambda\lambda'}$ is exactly the transformation rule
of the shear-bend coordinates of a pleated torus under diagonal
exchanges.

\begin{lem}[Lemma 27 in \cite{BonLiu}] \label{lem:nonquantum}
Let $q$ be a primitive $N$-th root
of unity with $N$ odd, and let $\rho:\mathcal
T^q_\lambda\to\mathrm{End}(V)$ be an irreducible representation
with
$$\rho(X_1^N)=x_1\,\mathrm{id}_V,\quad
\rho(X_2^N)=x_2\,\mathrm{id}_V, \quad
\rho(X_3^N)=x_3\,\mathrm{id}_V,\quad \rho(H)=h\,\mathrm{id}_V$$
then $\rho':=\rho\circ\Phi^q_{\lambda\lambda'}: \mathcal
T^q_{\lambda'}\to\mathrm{End}(V)$ is an irreducible representation
with $$\rho'\big(X_1'{}^N\big)=x_1'\,\mathrm{id}_V, \quad
\rho'\big(X_2'{}^N\big)=x_2'\,\mathrm{id}_V, \quad
\rho'\big(X_3'{}^N\big)=x_3'\,\mathrm{id}_V, \quad
\rho(H')=h'\,\mathrm{id}_V$$ where
$(x_1',x_2',x_3')=\Phi^1_{\lambda\lambda'}(x_1,x_2,x_3)$ and
$h'=h$. \qed
\end{lem}

Note: Since
$\Phi^q_{\lambda\lambda'}(\mathcal T^q_{\lambda'})\nsubseteq
\mathcal T^q_\lambda$, one needs  Lemma 25 and Lemma 26 in
\cite{BonLiu} to rigorously make sense of
$\rho\circ\Phi^q_{\lambda\lambda'}$ as a representation of $\mathcal
T^q_{\lambda'}$.
However, this turns out to coincide with our intuitive
understanding.

\section{Algebraic Lemmas}\label{sec:AlgLem}

Given a primitive $N$-th root of unity with $N$ odd, let $\mathcal
W_q$ be the algebra defined by generators $U^{\pm1}$, $V^{\pm1}$,
$W^{\pm1}$ and by the relations $VU=q^4UV$, $WV=q^4VW$, $UW=q^4WU$.
It is isomorphic to the Chekhov-Fock algebra $\mathcal T^q_\lambda$
for any ideal triangulation $\lambda$ of the 1-puncture torus. The
central element is $H=q^2UVW$. This will not cause any trouble
because we are going to deal with automorphisms of $\mathcal W$
instead of isomorphisms $\Phi^q_{\lambda\lambda'}$ between different
algebras.

\begin{lem}[{See for instance \cite[Sect.~4]{BonLiu}}]
\label{lem:std} Every finite-dimensional irreducible
representation of $\mathcal W_q$ is conjugate to $\chi_{u,v,h}:
\mathcal W_q\to\mathrm{End}(\mathbb C^N)$, defined by
\begin{equation*}
\chi_{u,v,h}(U)=u\left(
\begin{array}{ccccc}
1 & 0 & 0 & \cdots & 0 \\
0 & q^4 & 0 & \cdots & 0\\
0 & 0 & q^8 & \cdots & 0\\
\vdots & \vdots & \vdots & \ddots & \vdots\\
0 & 0 & 0 & \cdots & q^{4(N-1)}
\end{array}\right),
\quad \chi_{u,v,h}(V)=v\left(
\begin{array}{ccccc}
0 & 1 & 0 & \cdots & 0\\
0 & 0 & 1 & \cdots & 0\\
\vdots & \vdots & \vdots & \ddots & \vdots\\
0 & 0 & 0 & \cdots & 1\\
1 & 0 & 0 & \cdots & 0
\end{array}\right),\quad
\end{equation*}
and
$$ \chi_{u,v,h}(W)=\frac{q^{-2}h}{uv}\, \left(
\begin{array}{ccccc}
0 & 0 & 0 & \cdots & q^{-4(N-1)}\\
1 & 0 & 0 & \cdots & 0\\
0 & q^{-4} & 0 & \cdots & 0\\
\vdots & \vdots & \ddots & \vdots & \vdots\\
0 & 0 & \cdots & q^{-4(N-2)} & 0
\end{array}\right).$$
Note that in this case $\chi_{u,v,h}(H)=h\,\mathrm{id}_{\mathbb
C^N}$. \qed
\end{lem}

Let $\widehat{\mathcal W}_q$ be the fraction division algebra of
$\mathcal W_q$. Consider the automorphism $\mathcal R:
\widehat{\mathcal W}_q\to \widehat{\mathcal W}_q$ defined by
the property that
\begin{align} \label{eqn:R} \mathcal R(U)&=
\big(1+qU^{-1}\big)^{-1}\big(1+q^3U^{-1}\big)^{-1}W\nonumber\\
\mathcal R(V)&= \big(1+qU\big)\big(1+q^3U\big)V\\
\mathcal R(W)&= U^{-1}.\nonumber
\end{align}
One easily see that $\mathcal R$ does define an algebra
isomorphism, namely that $\mathcal R(U)$, $\mathcal R(V)$ and
$\mathcal R(W)$ satisfy the same relations as $U$, $V$ and $W$
(compare Proposition 5 in \cite{Liu}), and that $\mathcal R$ is
invertible.

\begin{lem} \label{lem:R}
If $\chi_{u,v,h}$ is a standard representation as in
Lemma~\ref{lem:std}, then
$$\chi_{u,v,h}\circ\mathcal R(X)= C_R\cdot\chi_{u',v',h'}(X)\cdot
C_R^{-1}, \qquad \forall\ X\in\mathcal W_q,$$ where
$$(u')^N=\frac{h^N}{u^Nv^N(1+u^{-N})^2},
\qquad (v')^N=(1+u^N)^2\,v^N,\qquad h'=h$$ and where the matrix
$C_R(u,v,u',v',h)\in\mathrm{GL}_N(\mathbb C)$ is defined by its
entries
$$ (C_R)_{ij}=q^{2(j-i)^2}\
\Big(\frac{u'v'}{uh}\Big)^{j-i}\
\Big(\frac{v'}{v}\Big)^i\
\prod_{\alpha=1}^{i}\frac{1}{\big(1+q^{4\alpha-3}u\big)
\big(1+q^{4\alpha-1}u\big)}.$$
\end{lem}

\begin{proof} Lemma~\ref{lem:nonquantum} shows that
$\chi_{u,v,h}\circ\mathcal R$ is isomorphic to a standard
representation $\chi_{u',v',h'}$ with
$(u')^N=\frac{h^N}{u^Nv^N(1+u^{-N})^2}$,
$(v')^N=(1+u^N)^2\,v^N$ and $h'=h$.
To compute the conjugation matrix $C_R$, the last two
equalities of (\ref{eqn:R}) provide
\begin{align*}
(C_R)_{ij}&=\frac{v'}{v}\,\frac{1}{(1+q^{4i-3}u)(1+q^{4i-1}u)}\,
(C_R)_{i-1\ j-1}\\
(C_R)_{ij}&=q^{-4i}\ q^{4(j-1)}\ \frac{u'v'}{uq^{-2}h}\ (C_R)_{i\
j-1}
\end{align*}
and these inductive relations immediately give the entries of
$C_R$.
\end{proof}

Similarly, let the isomorphism $\mathcal L:\widehat{\mathcal
W}_q\to
\widehat{\mathcal W}_q$ be defined by
\begin{align} \label{eqn:L}
\mathcal L(U)&=
\big(1+qV^{-1}\big)^{-1}\big(1+q^3V^{-1}\big)^{-1}U\nonumber\\
\mathcal L(V)&= \big(1+qV\big)\big(1+q^3V\big)W\\
\mathcal L(W)&= V^{-1},\nonumber
\end{align}

\begin{lem} \label{lem:L}
If $\chi_{u,v,h}$ is a standard representation as in
Lemma~\ref{lem:std}, then
$$\chi_{u,v,h}\circ\mathcal L(X)= C_L\cdot\chi_{u'',v'',h''}(X)\cdot
C_L^{-1}, \qquad \forall\ X\in\mathcal W_q,$$ where
$$(u'')^N=\frac{u^N}{(1+v^{-N})^2},\qquad (v'')^N=
\frac{h^N(1+v^N)^2}{u^Nv^N},\qquad h''=h$$
where
$C_L=G\widetilde{C}_L$ is the product of matrices
$G$, $\widetilde{C}_L(u,v,u'',v'',h)\in\mathrm{GL}_N(\mathbb C)$
with entries $G_{ij}=q^{4ij}$
and
$$ (\widetilde{C}_L)_{ij}=q^{2(j-i)^2+2i^2}\ \Big(\frac{u''v''}{vh}\Big)^{j-i}\
\Big(\frac{uvv''}{h}\Big)^i\
\prod_{\alpha=1}^{i}\frac{1}{\big(1+q^{4\alpha-3}v\big)
\big(1+q^{4\alpha-1}v\big)}.           $$
\end{lem}

\begin{proof}
The proof is similar to the proof of the previous lemma, except
that $\chi_{u,v,h}(V)$ is not diagonal. This inconvenience is
bypassed in the following way. Define another irreducible
representation $\mu_{u,v,h}$ by
\begin{equation*}
\mu_{u,v,h}(U)=u \left(
\begin{array}{ccccc}
0 & 0 & 0 & \cdots & 1\\
1 & 0 & 0 & \cdots & 0\\
0 & 1 & 0 & \cdots & 0\\
\vdots & \vdots & \ddots & \vdots & \vdots\\
0 & 0 & \cdots & 1 & 0
\end{array}\right), \quad
\mu_{u,v,h}(V)=v\left(
\begin{array}{ccccc}
1 & 0 & 0 & \cdots & 0 \\
0 & q^4 & 0 & \cdots & 0\\
0 & 0 & q^8 & \cdots & 0\\
\vdots & \vdots & \vdots & \ddots & \vdots\\
0 & 0 & 0 & \cdots & q^{4(N-1)}
\end{array}\right),
\end{equation*}
\begin{equation*}
\mu_{u,v,h}(W)=\frac{q^{-2}h}{uv}\,\left(
\begin{array}{ccccc}
0 & 1 & 0 & \cdots & 0\\
0 & 0 & q^{-4} & \cdots & 0\\
\vdots & \vdots & \vdots & \ddots & \vdots\\
0 & 0 & 0 & \cdots & q^{-4(N-2)}\\
q^{-4(N-1)} & 0 & 0 & \cdots & 0
\end{array}\right).\quad
\end{equation*}
One easily computes that $G\cdot\mu_{u,v,h}(X)\cdot
G^{-1}=\chi_{u,v,h}(X)$, with $G_{ij}=q^{4ij}$. Let
$\widetilde{C}_L$ be the matrix such that
$\mu_{u,v,h}\circ\mathcal L(X)=\widetilde{C}_L\cdot
\chi_{u'',v'',h}(X)\cdot\widetilde{C}_L^{-1}$. Obviously
$C_L=G\widetilde{C}_L$. As in the previous lemma, we obtain
the inductive equations
\begin{align*}
(\widetilde{C}_L)_{ij}&=q^{4(i-1)}\ \Big(\frac{uvv''}{q^{-2}h}\Big)\
\frac{1}{\big(1+q^{4i-3}v\big)\big(1+q^{4i-1}v\big)}\
(\widetilde{C}_L)_{i-1\ j-1}\\
(\widetilde{C}_L)_{ij}&=q^{-4i}\ q^{4(j-1)}\
\frac{u''v''}{vq^{-2}h}\ (\widetilde{C}_L)_{i\ j-1}
\end{align*}
which provide the entries of $\widetilde{C}_L$.
\end{proof}

\begin{rem}
  Note that, although the  quantities $u$, $v$, $h$,
$u'$,
$v'$,  involved in the definition of the matrix
$C_R$ are related by equations involving only their $N$--th
powers, the matrix itself depends on the choice of these
$N$--roots. The same holds for $C_L$. In the next
section, the $N$--powers will be determined by geometric
information but the actual computation will depend on actual
choices of
$N$--roots for this geometric data.
\end{rem}

\section{The invariant in the case of the $1$--puncture torus}
\label{sec:1-PuncTorus}

We now compute the invariant $C_\varphi$ of
Theorem~\ref{thm:invariant} in the case of the 1-puncture torus $T$.
The computation is here greatly simplified by the fact that the
mapping class group of $T$ is isomorphic to $\mathrm{SL}_2(\mathbb
Z)$. In addition, recall that a mapping class is pseudo-Anosov
exactly when the corresponding element of $\mathrm{SL}_2(\mathbb Z)$
has a trace of absolute value greater than 2.

\begin{prop} \label{prop:decomp}
Let $A$ be an element of $\mathrm{SL}_2(\mathbb{Z})$ with
$|\mathrm{Tr}(A)|>2$. Then the conjugacy class of $A$ in
$\mathrm{SL}_2(\mathbb Z)$ contains an element of the form
$$
\left(\begin{array}{cc}1 & a_1\\0 & 1\end{array}\right)
\left(\begin{array}{cc}1 & 0\\b_1 & 1\end{array}\right)
\left(\begin{array}{cc}1 & a_2\\0 & 1\end{array}\right)
\left(\begin{array}{cc}1 & 0\\b_2 & 1\end{array}\right)\ \cdots\
\left(\begin{array}{cc}1 & a_n\\0 & 1\end{array}\right)
\left(\begin{array}{cc}1 & 0\\b_n & 1\end{array}\right) $$ where
$n>0$ and the $a_i$, $b_i$ are positive integers. Moreover, the
right hand side is unique up to cyclic permutation of the factors
$\left(\begin{array}{cc}1 & a_i\\0 & 1\end{array}\right)$ and
$\left(\begin{array}{cc}1 & 0\\b_i & 1\end{array}\right)$. \qed
\end{prop}

In other words, after conjugation, every mapping class of $T$ can be
represented by a matrix $A = A_1A_2\cdots A_n$ with $A_i=R$ or $L$ ,
where
  $L= \left(\begin{array}{cc}1 & 1\\0 & 1\end{array}\right)$ and
$R= \left(\begin{array}{cc}1 & 0\\1 & 1\end{array}\right)$.

The identification between $H_1(T) \cong \mathbb Z^2$ assigns a
slope in $\mathbb Q\cup\{\infty\}$ to each edge of an ideal
triangulation $\lambda$ of $T$.
Consider the ideal triangulation $\lambda_{(0)}$ whose edges have
respective slopes 0,
$\infty$ and 1. Then define ideal triangulations
$\lambda_{(0)}$, $\lambda_{(1)}$, \dots,
$\lambda_{(n)}=\varphi(\lambda_{(0)})$ by the property that
$\lambda_{(i)} = A_1A_2 \dots A_i(\lambda_{(0)})$.

It is
immediate that the edges $\lambda_{(i)1}$, $\lambda_{(i)2}$ and
$\lambda_{(i)3}$ of $\lambda_{(i)}$ all have non-negative slope.
Choose the indexing of these edges so that $\lambda_{(i)1}$ has
the lowest slope and $\lambda_{(i)2}$ has the highest slope of
the three. Note that with this convention the edges
$\lambda_{(i)1}$,
$\lambda_{(i)2}$ and
$\lambda_{(i)3}$ occur counterclockwise in this order around the
two triangles of $\lambda_{(i)}$.

The key observation is the following:
\begin{itemize}
\item if $A_i=R$, then $\lambda_{(i)}$ is obtained from
$\lambda_{(i-1)}$ by a diagonal exchange along the edge
$\lambda_{(i-1)1}$, followed by a reindexing exchanging
$\lambda_{(i)1}$ and $\lambda_{(i)3}$;
\item if $A_i=L$, then $\lambda_{(i)}$ is obtained from
$\lambda_{(i-1)}$ by a diagonal exchange along the edge
$\lambda_{(i-1)2}$, followed by a reindexing exchanging
$\lambda_{(i)2}$ and $\lambda_{(i)3}$.
\end{itemize}

\begin{figure}[h]
\begin{center}
\includegraphics{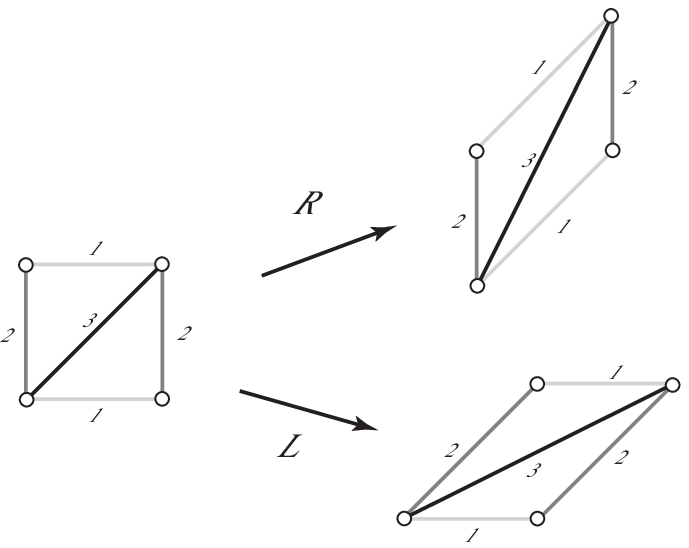}
\end{center}
\caption{} \label{fig:LR}
\end{figure}

This property is illustrated in
Figure~\ref{fig:LR} for the case where $i=1$, and the other
cases follow from this one by observing that
$$
\lambda_{(i)} = (A_1A_2 \dots A_{i-1}) A_i (A_1A_2 \dots
A_{i-1})^{-1} (\lambda_{(i-1)}).
$$
Incidentally, this accounts for the somewhat unnatural order of
the $A_k$ in the definition of $\lambda_{(i)}$.

For the complete hyperbolic metric on the mapping torus
$M_{\varphi}=S\times[0,1]/\sim
$ with its complete hyperbolic metric, each ideal
triangulation
$\lambda_{(i)}$ determines a unique pleated surface
$f_{\lambda_{(i)}}: T\to M_\varphi$ with
pleating locus $\lambda_{(i)}$.
The geometry of this pleated surface associates complex weights
$x_{(i)1}$, $x_{(i)2}$, $x_{(i)3}\in \mathbb C-\{0\}$ to the
edges of $\lambda_{(i)}$, corresponding to the exponential
shear-bend coordinates of $f_{\lambda_{(i)}}$ along the
components of its pleating locus.

\begin{rem}
  In \cite{Gue}, Gu\'eritaud proves the remarkable fact that the
pleated surfaces $f_{\lambda_{(i)}}$ are embeddings. In
addition, one passes from$f_{\lambda_{(i-1)}}$ to
$f_{\lambda_{(i)}}$ by a diagonal exchange across a
positively oriented hyperbolic ideal tetrahedra with dihedral
angles strictly between $0$ and
$\pi$.
\end{rem}

\begin{thm}
Given any pseudo-Anosov diffeomorphism $\varphi: T\to T$, if
$\varphi$ has the canonical decomposition $A_1A_2\cdots A_n$ with
$A_i=R$ or $L$ as in Proposition~\ref{prop:decomp}, then the
conjugacy class of $C_\varphi$ has a matrix representative
$$C_\varphi=C_1 C_2\cdots C_n $$ where

\begin{itemize}

\item if $A_i=R$, then $C_i=C_R(u_{i-1},v_{i-1},u_i,v_i,1)$, defined
in Lemma~\ref{lem:R}, with
$$(u_i)^N=\frac{1}{u_{i-1}^Nv_{i-1}^N(1+u_{i-1}^{-N})^2}, \qquad
(v_i)^N=(1+u_{i-1}^N)^2\,v_{i-1}^N,$$

\item if $A_i=L$, then $ C_i=C_L(u_{i-1},v_{i-1},u_i,v_i,1)$,
defined in Lemma~\ref{lem:L}, with
$$(u_i)^N=\frac{u_{i-1}^N}{(1+v_{i-1}^{-N})^2},\qquad (v_i)^N=
\frac{(1+v_{i-1}^N)^2}{u_{i-1}^Nv_{i-1}^N},$$

\end{itemize}

and, the initial values $u_0$, $v_0$, $h$ are chosen so that
$$ u_n=u_0, \qquad v_n=v_0, \qquad h=1.$$

\end{thm}

\begin{proof}
The action of the word $A_1A_2\cdots A_n$ produces a sequence of
ideal triangulations
$$ \lambda=\lambda_{(0)}\stackrel{A_1}{\longrightarrow} \lambda_{(1)}
\stackrel{A_2}{\longrightarrow} \lambda_{(2)}\longrightarrow
\cdots\stackrel{A_n}{\longrightarrow}\lambda_{(n)}=
\varphi(\lambda),$$ which in turn gives a composition of
isomorphisms

$$ \xymatrix{ \widehat{\mathcal T}^q_{\varphi(\lambda)} &
\ar[rr]^{\Phi^q_{\lambda_{(n-1)}\lambda_{(n)}}} & & &
\widehat{\mathcal T}^q_{\lambda_{(n-1)}}\ \longrightarrow\ \cdots\ &
\ar[rr]^{\Phi^q_{\lambda_{(0)}\lambda_{(1)}}} & & &
\widehat{\mathcal T}^q_\lambda}.$$

By sending $X_{\lambda_{(i)1}}$, $X_{\lambda_{(i)2}}$,
$X_{\lambda_{(i)3}}$ respectively to $U$, $V$, $W$ we identify each
$\widehat{\mathcal T}^q_{\lambda_{(i)}}$ to the algebra
$\widehat{\mathcal W}_q$ of \S\ref{sec:AlgLem}. Then we get an
automorphism $\mathcal A_1\mathcal A_2\cdots \mathcal
A_n:\widehat{\mathcal W}_q \to \widehat{\mathcal W}_q$ with
$\mathcal A_i\in\{\mathcal R,\mathcal L\}$. This automorphism is
conjugate to the automorphism
$\Phi^{-1}\circ\Phi^q_{\lambda\varphi(\lambda)}:\widehat{\mathcal
T}^q_{\varphi(\lambda)}\to \widehat{\mathcal T}^q_{\varphi(\lambda)}
$, ($\Phi:\widehat{\mathcal T}^q_{\varphi(\lambda)} \to
\widehat{\mathcal T}^q_\lambda$ is the natural identification). If
we choose a standard representation $\chi_{u,v,h}$, then
Lemma~\ref{lem:R} and Lemma~\ref{lem:L} tell us that for any $i$,
$1\le i\le n$, we have
\begin{align*}
{} &\chi_{u,v,h}\circ \mathcal A_1\mathcal A_2 \cdots\mathcal
A_i(X)\\
{}=& C_1(u_0,v_0,u_1,v_1,h)\cdot \Big(\chi_{u_1,v_1,h}\circ \mathcal
A_2\cdots
\mathcal A_i(X)\Big)\cdot C_1(u_0,v_0,u_1,v_1,h)^{-1}\\
{}=&\cdots\\
{}=& (C_1\cdots C_i)\Big(\chi_{u_k,v_k,h}(X)\Big) (C_1\cdots
C_i)^{-1}
\end{align*} for all $X\in\mathcal W_q$.
In particular, the hyperbolic metric of the mapping torus
$M_\varphi$ provides a interesting initial standard representation
$\chi_{u_0,v_0,h_0}$ such that $u_n^N=u_0^N$, $v_n^N=v_0^N$ (We have
interpreted $u_i^N$, $v_i^N$ to be the shear-bend coordinates
$x_{(i)1}$, $x_{(i)2}$). We choose the standard representations from
their conjugacy classes so that $u_n=u_0$, $v_n=v_0$, and we choose
$h=1$ because it is natural to work on the ``cusped quantum
Teichm\"uller space'' (see \cite{BonLiu, Liu}) in which we should
put the ``cusp condition'' $H^2=1$, and a geometric consideration
fixes $H=1$. The above equation applied to $i=n$ gives
\begin{align*}
{}&\chi_{u_0,v_0,h_0}\circ \mathcal A_1\mathcal A_2 \cdots\mathcal
A_n(X)\\
{}=& (C_1\cdots C_n)\Big(\chi_{u_0,v_0,h_0}(X)\Big) (C_1\cdots
C_n)^{-1}
\end{align*}
Comparing with equation~\ref{eqn:invariant}, the product of matrices
$C_1\cdots C_n$ turns out to be a representative of the invariant
$C_\varphi$.

\end{proof}

\section{The invariant in the case of the 4-puncture sphere}
\label{sec:4-PuncSphere}

We identify the 4-puncture sphere to be the quotient $S\cong
(\mathbb R^2-\mathbb Z^2)/G$ where $G$ is the group generated by
four elements $ (x,y)\mapsto (-x,y)$, $ (x,y)\mapsto (x+2,y)$,
$(x,y)\mapsto (x,-y)$ and 
$(x,y)\mapsto (x,y+2)$. The mapping class group is the semi-directed
product of $\mathrm{PSL}_2(\mathbb Z)$ on $\mathbb Z/2\mathbb Z\
\oplus\ \mathbb Z/2\mathbb Z$. The image of a diffeomorphism
$\varphi$ can be represented by a matrix
$A_\varphi\in\mathrm{SL}_2(\mathbb Z)$.

\begin{figure}[h]
\begin{center}
\includegraphics{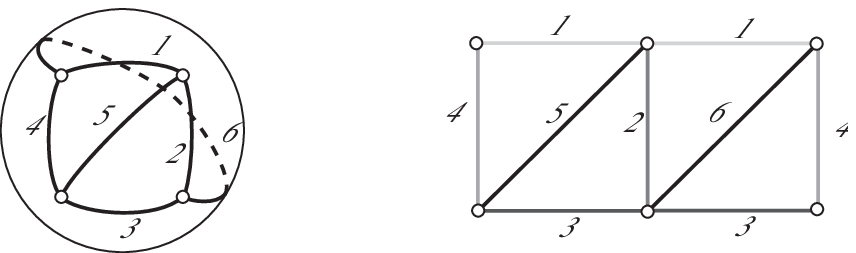}
\caption{An ideal triangulation $\lambda$ of $S$.}
\end{center}
\label{fig:sphere}
\end{figure}

We will restrict attention to ideal triangulation of $S$ which are
isomorphic to the 2-skeleton of an ideal tetrahedron, as in the
figure above. This is equivalent to the property that each puncture
is adjacent to exactly 3 edges.

If $\varphi: S\to S$ is associated to the matrix
$A_\varphi\in\mathrm SL_2(\mathbb Z)$, then $\varphi$ is
pseudo-Anosov exactly when $|\mathrm{Tr}\, \varphi|>2$. Similar to
the 1-puncture torus, such a matrix has an $LR$ decomposition. The
difference now is that, to stay within tetrahedral ideal
triangulations, the $L$ and $R$ are now associated to the product of
two diagonal exchanges along opposite edges. More precisely, with
the indexing convention of the figure above, $L$ corresponds to
$\Delta_2\Delta_4$ and $R$ corresponds to $\Delta_1\Delta_3$. By
\cite{Liu}, the Chekhov-Fock algebra $\mathcal T^q_\lambda$
($\lambda$ as in the figure) is generated by $X_1$, $X_2$ and the
central elements
$$ P_1=qX_1X_2X_5,\qquad P_2=q^{-1}X_2X_3X_6,
\qquad P_3=qX_3X_4X_5,\qquad P_4=qX_1X_4X_6, $$
$$ H=q^2X_1X_2X_3X_4X_5X_6.$$
There exist automorphisms $\mathcal R$ and $\mathcal L$ of the
Chekhov-Fock algebra $\mathcal T^q_\lambda$ which fix the central
elements and map $X_1$, $X_2$ as follows,
\begin{align*}
\mathcal R(X_1)&=\big(1+qX_1^{-1}\big)^{-1}\big(1+qX_3^{-1}\big)^{-1}X_6,\\
\mathcal R(X_2)&=(1+qX_1)(1+qX_3)X_2,
\end{align*}
and
\begin{align*}
\mathcal L(X_1)&=\big(1+qX_2^{-1}\big)^{-1}\big(1+qX_4^{-1}\big)^{-1}X_1,\\
\mathcal L(X_2)&=(1+qX_2)(1+qX_4)X_5.
\end{align*}

If we use the standard representation
\begin{equation*}
\chi_{u,v,h,p}(X_1)=u\left(
\begin{array}{ccccc}
1 & 0 & 0 & \cdots & 0 \\
0 & q^2 & 0 & \cdots & 0\\
0 & 0 & q^4 & \cdots & 0\\
\vdots & \vdots & \vdots & \ddots & \vdots\\
0 & 0 & 0 & \cdots & q^{2(N-1)}
\end{array}\right),
\quad \chi_{u,v,h,p}(X_2)=v\left(
\begin{array}{ccccc}
0 & 1 & 0 & \cdots & 0\\
0 & 0 & 1 & \cdots & 0\\
\vdots & \vdots & \vdots & \ddots & \vdots\\
0 & 0 & 0 & \cdots & 1\\
1 & 0 & 0 & \cdots & 0
\end{array}\right),\quad
\end{equation*}
and
\begin{equation*}
\chi_{u,v,h,p}(H)= h\ \mathrm{id}_{\mathbb C^N},\quad
\chi_{u,v,h,p}(P_j)= p_j\ \mathrm{id}_{\mathbb C^N},\quad
j=1,2,3,4,
\end{equation*}
then we have similar algebraic Lemmas for the 4-puncture sphere,
with the new matrix functions,
\begin{align*}
(C^*_R)_{ij} & = q^{(i-j)(i-j-2)} \left(\frac{v'}{v}\right)^j 
(uvp_3u')^{j-i}\ \prod_{\alpha=1}^{i-j}\frac{1}{(1+q^{1-2\alpha} u^{-1})
\left(1+\frac{q^{-1-2\alpha}h}{up_2p_3}\right)}\quad \times \\
 & \qquad \times\quad \prod_{\beta=i-j+1}^i 
\frac{1}{(1+q^{2\beta-1}u)
\left(1+ \frac{q^{2\beta+1}up_2p_3}{h}\right)} \\
& {} \\
G^*_{ij} & =\ q^{2ij} \\
& {} \\
(\widetilde{C}^*_L)_{ij} & =q^{j^2}\left(\frac{u'}{u}\right)^{j-i}
\left(\frac{uvv'}{p_1}\right)^j\qquad \prod_{\alpha=1}^{i-j}
\frac{1}{(1+q^{1-2\alpha}v^{-1})\left(1+\frac{q^{1-2\alpha}h}
{vp_3p_4}\right)}\quad \times \\
& \qquad \times\quad\prod_{\beta=i-j+1}^{i}\frac{1}{(1+q^{2\beta-1}v)
\left(1+\frac{q^{2\beta-1}vp_3p_4}{h}\right)}
\end{align*}

\vspace{0.5cm}

The same arguments as in \S\ref{sec:1-PuncTorus} then give:

\begin{thm}

Given any pseudo-Anosov diffeomorphism $\varphi: S\to S$, if
$\varphi$ has the canonical decomposition $A_1A_2\cdots A_n$ with
$A_i=R$ or $L$ as in Proposition~\ref{prop:decomp}, then the
conjugacy class of $C_\varphi$ has a matrix representative
$$C^*_\varphi=C^*_1 C^*_2\cdots C^*_n $$ where

\begin{itemize}

\item if $A_i=R$, then $C^*_i=C^*_R(u_{i-1},v_{i-1},u_i,v_i)$, defined
above, with
$$(u_i)^N=-\,\frac{1}{u_{i-1}^Nv_{i-1}^N(1+u_{i-1}^{-N})^2}, \qquad
(v_i)^N=(1+u_{i-1}^N)^2\,v_{i-1}^N,$$

\item if $A_i=L$, then $ C^*_i=C^*_L(u_{i-1},v_{i-1},u_i,v_i,1)=G^*\cdot
\widetilde{C}^*_L(u_{i-1},v_{i-1},u_i,v_i,1)$, defined above, with
$$(u_i)^N=\frac{u_{i-1}^N}{(1+v_{i-1}^{-N})^2},\qquad (v_i)^N=
-\,\frac{(1+v_{i-1}^N)^2}{u_{i-1}^Nv_{i-1}^N},$$

\end{itemize}

and, the initial values $u_0$, $v_0$, $h$, $p$ are chosen so that
$$ u_n=u_0, \qquad v_n=v_0, \qquad h=1, \qquad p_j=1,\ \forall\ j.$$
\qed

\end{thm}

\end{document}